# Factoring Catalan numbers

Gennady Eremin
ergenns@gmail.com

August 11, 2019

**Abstract.** The paper describes a prime factorization of the Catalan numbers. Odd prime factors are distributed in layers in accordance with Legendre's formula. The content of each layer is a network of the intervals, Chebyshev's Segments. The primes of Segment are not calculated and are selected on the basis of its bounds. Layers contain non-repeated primes. Repeated factors are formed when primes are duplicated among different layers. The paper slightly modifies Kummer's theorem for the selection of individual prime factors, also starting from the boundaries of Segments. In conclusion, the reader is offered a software service for factorization of the Catalan number with index up to $10^8$.

**Keywords.** Catalan number, prime factorization, Legendre's formula, Chebyshev's Segments, Kummer's theorem.

## 1 Introduction

Currently, the decomposition of natural numbers is relevant and in demand to work with big numbers, such as in cryptography. In cryptanalysis, we factorize keys whose length is hundreds of decimal digits. In this regard, the decomposition of special numbers is interesting, since mathematical formulas (recurrent relations, etc) significantly prime factorization of big integers of known numerical sequences. To factorize a huge special number, you may not know its natural form. For example, the 9,999-th Catalan number has a length of 6,000 decimal digits. Of course, such a number is almost impossible to get, but it is easy to factorize. In this article, we will factorize the Catalan numbers with large indexes in examples.

**1.1**. The Catalan numbers appear in many combinatorial applications [St15]. Here is the known explicit formula for the general term of the Catalan sequence:

$$(1.1) \qquad Cat(n) \;=\; \binom{2n}{n}\Big/(n+1) \;=\; \frac{(2n)!}{n!(n+1)!}, \quad n \geq 0.$$

The first Catalan numbers for $n = 0, 1, 2, 3 \ldots$ are (see [A108])

  1, 1, 2, 5, 14, 42, 132, 429, 1430, 4862, 16796, 58786, 208012…

In number theory for a prime $p$, the *p-adic valuation* (or *order*) of $n \in \mathbb{N}$, $v_p(n)$, is the highest power of $p$ that divides $n$ (or the number of factors of $p$ in the prime factorization of $n$). For example, $v_2(6) = v_2(2) = 1$, $v_3(36) = v_3(9) = 2$, $v_5(1000) = v_5(125) = 3$. Let's write some obvious properties for a prime $p$ and $n, m, k \in \mathbb{N}$:

$$v_p(nm) = v_p(n) + v_p(m), \quad v_p(n/m) = v_p(n) - v_p(m), \quad v_p(p^k) = k.$$

It is easy to convert (1.1) to a more convenient form [We16]

(1.2) $$Cat(n) = 2^n \times (2n-1)!! / (n+1)!,$$

where $(2n-1)!! = 1 \cdot 3 \cdot 5 \cdots (2n-1)$ is a double odd factorial. From (1.2) we get the 2-adic valuation of the $n$-th Catalan number: $v_2(Cat(n)) = n - v_2((n+1)!)$.

The latter equality can be further simplified (see Theorem 2.1 in [DS18])

(1.3) $$v_2(Cat(n)) = wt(n+1) - 1,$$

where $wt(n+1)$ is the sum (*weight*) of the digits in the base 2 expansion of $n+1$, i.e., the sum of the units of the binary code of $n+1$. Obviously, the $n$-th Catalan numbers are odd if $wt(n+1) = 1$, or $n = 2^k - 1, k \in \mathbb{N}$, i.e. $n = 1, 3, 7, 15$, etc.

**Example 1.1.** Start factoring the 9,999-th Catalan number, in the decomposition of which there are 1,560 primes including powers (trust us, reader). First, check an even prime 2. Let's calculate the 2-adic valuation of $Cat(9,999)$. Below we convert the decimal integer 10,000 to binary word with 5 units.

$$v_2(Cat(9,999)) = wt(9,999 + 1) - 1 = wt(10,011,100,010,000_2) - 1 = 5 - 1 = 4.$$

It remains to choose odd prime factors in the amount of $1560 - 4 = 1556$. □

**1.2.** From (1.1) directly follows, first, prime factors of $Cat(n)$ are less than $2n$. Let us say that we are dealing with the finite *prime interval* $_p(1; 2n) = \{2, 3, 5 \ldots\}$. We call this interval the *factor base* of the $n$-th Catalan number. Any prime interval, *P-interval*, contains only primes within the specified boundaries. One more thing, any $p \in {_p(n+1; 2n)}$ divides $Cat(n)$, but $p^2 \nmid Cat(n)$. That is, primes from $_p(n+1; 2n)$ are non-repeating factors of the $n$-th Catalan number.

**Note 1.2.** In general, a prime interval can have both open and closed boundaries. Boundaries can be arbitrary positive real numbers. For example, $_p[3; 8] = \{3, 5, 7\}$, $_p(4.99; 11] = \{5, 7, 11\}$, $_p(\sqrt{3}; 11.001) = \{2, 3, 5, 7, 11\}$, $_p[13.01; 17) = \emptyset$.

**Example 1.3.** Continue factoring the 9,999-th Catalan number. The prime interval $_p(10000; 19998) = \{10007, 10009, \ldots, 19997\}$ contents 1,033 non-repeating factors of $Cat(9,999)$. We still have to select $1556 - 1033 = 523$ odd prime factors, which are distributed over other prime intervals. □

The $n$-th Catalan number and the binomial coefficient $\binom{2n}{n}$ have the same composition of prime factors. In 1850 Pafnuty Chebyshev noted that $\binom{2n}{n}$ is divisible by the product of the primes in the interval $(n, 2n)$ [Po15]. Therefore, let's call the interval $_p(n+1; 2n)$ *Chebyshev's Segment*. The first results of the work with Chebyshev's Segments are described in [Er16].

# 2 Chebyshev's Segments and Black Holes

**2.1.** The first (main) Chebyshev's Segment, $S_1(n) = {}_p(n+1; 2n)$, contains most of the prime factors of the $n$-th Catalan number. The remaining primes are distributed to other Segments between which inaccessible gaps, *Black Holes*, are formed. Chebyshev's Segments do not intersect and have open boundaries similar to $S_1(n)$. Black Holes also do not intersect, and it is logical to use closed boundaries for them. In this case, the value of the boundaries for adjacent Segments and Black Holes will be the same. For example, the first Black Hole $H_1(n) = {}_p[?; n+1]$ is adjacent tightly to $S_1(n)$. Segments and Black Holes are often empty. For example, $S_1(2) = {}_p(3; 4) = \emptyset$, but $S_1(3) = {}_p(4; 6) = \{5\}$. Let's describe (1.2) in more detail:

(2.1) $\quad Cat(n) = 2^n \times A/B, \ A = 1 \cdot 3 \cdot 5 \cdots (2n-1), \ B = 1 \cdot 2 \cdot 3 \cdots (n+1).$

For an odd prime $p$,

$$v_p(Cat(n)) = v_p(A/B) = v_p(A) - v_p(B).$$

Let $A_1$ denote P-interval in which each prime number $p$ occurs in A only once, i.e. $v_p(A) = 1$. Obviously, $A_1 = {}_p(2n/3; 2n)$. Accordingly for B, there is the similar P-interval $B_1 = {}_p((n+1)/2; n+1]$; for each $p \in B_1$, $v_p(B) = 1$. The sets $A_1$ and $B_1$ intersect and we can get the first Chebyshev's Segment and the first Black Hole:

$S_1(n) = A_1 \setminus B_1 = {}_p(2n/3; 2n) \setminus {}_p((n+1)/2; n+1] = {}_p(n+1; 2n),$
$H_1(n) = A_1 \cap B_1 = {}_p(2n/3; 2n) \cap {}_p((n+1)/2; n+1] = {}_p[2n/3; n+1].$

Similarly, we obtain the second Chebyshev's Segment and the second Black Hole. Denote by $A_2$ the P-interval in which each prime $p$ occurs in A exactly twice, i.e. $v_p(A) = 2$. It is easy to see, $A_2 = {}_p(2n/5; 2n/3)$. The appropriate P-interval for B is $B_2 = {}_p((n+1)/3; (n+1)/2]$; for each $p \in B_2$, $v_p(B) = 2$. We calculate the second Chebyshev's Segment a bit differently.

$S_2(n) = A_2 \cap B_1 = {}_p(2n/5; 2n/3) \cap {}_p((n+1)/2; n+1] = {}_p((n+1)/2; 2n/3),$
$H_2(n) = A_2 \cap B_2 = {}_p(2n/5; 2n/3) \cap {}_p((n+1)/3; (n+1)/2] = {}_p[2n/5; (n+1)/2].$

It is logical to consider $B_0 = {}_p(n+1; \infty]$. On the basis of the last equations we give general formulas for Chebyshev's Segment and Black Hole:

$S_t(n) = A_t \cap B_{t-1}$
$\quad = {}_p(2n/(2t+1)); \ 2n/(2t-1)) \cap {}_p(((n+1)/t; \ (n+1)/(t-1)]$
$\quad = {}_p((n+1)/t; \ 2n/(2t-1)),$

$$H_t(n) = A_t \cap B_t$$
$$= {}_p(2n/(2t+1)); 2n/(2t-1)) \cap {}_p(((n+1)/(t+1); (n+1)/t]$$
$$= {}_p[2n/(2t+1); (n+1)/t], \quad t \geq 1.$$

As a result, we get

(2.2a) $$S_t(n) = {}_p((n+1)/t; 2n/(2t-1)), \quad t \in \mathbb{N}.$$

(2.2b) $$H_t(n) = {}_p[2n/(2t+1); (n+1)/t], \quad t \geq 0.$$

**2.2.** From (2.2) follows $H_t(n) \cup S_t(n) = A_t$. This means that Chebyshev's Segments and Black Holes of the form (2.2) cut the range ${}_p(2; 2n)$. Recall that we work with odd primes (the multiplicity for 2 we got in the previous section). The question naturally arises, how many Chebyshev's Segments are there? With the growth of $t$ (approaching the origin) Segments are sharply reduced, and in the limit for Segment boundaries must collapse (and then Segments are inverted). Let's equate the boundaries for Chebyshev's Segment in (2.2a):

$$(n+1)/t = 2n/(2t-1), \text{ then } t = (n+1)/2.$$

In this case for odd $n$, we get an empty Segment with the same boundaries ${}_p(2; 2)$. And the next Segment is inverted (the lower boundary exceeds the upper one). So we can formulate

**Proposition 2.1.** *For the n-th Catalan number, the number of Segments of type (2.2a) is less than $(n+1)/2$.*

We will omit the index $n$, if we know what the Catalan number we are talking about.

**Example 2.2.** Continue Example 1.3. Let's show the first 15 Segments for the 9,999-th Catalan number (the hash sign '#' denotes cardinality).

$S_1 = {}_p(10000; 19998) = \{10007, 10009, \ldots, 19997\}$, $\#S_1 = 1033$;
$S_2 = {}_p(5000; 6666) = \{5003, 5009, \ldots, 6661\}$, $\#S_2 = 190$;
$S_3 = {}_p(3333.3; 3999.6) = \{3343, 3347, \ldots, 3989\}$, $\#S_3 = 80$;
$S_4 = {}_p(2500; 2856.8) = \{2503, 2521, \ldots, 2851\}$, $\#S_4 = 47$;
$S_5 = {}_p(2000; 2222) = \{2003, 2011, \ldots, 2221\}$, $\#S_5 = 28$;
$S_6 = {}_p(1666.6; 1818) = \{1667, 1669, \ldots, 1811\}$, $\#S_6 = 19$;
$S_7 = {}_p(1428.5; 1538.3) = \{1429, 1433, \ldots, 1531\}$, $\#S_7 = 17$;
$S_8 = {}_p(1250; 1333.2) = \{1259, 1277, \ldots, 1327\}$, $\#S_8 = 13$;
$S_9 = {}_p(1111.1; 1176.3) = \{1117, 1123, \ldots, 1171\}$, $\#S_9 = 7$;
$S_{10} = {}_p(1000; 1052.5) = \{1009, 1013, \ldots, 1051\}$, $\#S_{10} = 9$;
$S_{11} = {}_p(909.09; 952.2) = \{911, 919, \ldots, 947\}$, $\#S_{11} = 6$;
$S_{12} = {}_p(833.3; 869.4) = \{839, 853, 857, 859, 863\}$, $\#S_{12} = 5$;
$S_{13} = {}_p(769.2; 799.9) = \{773, 787, 797\}$, $\#S_{13} = 3$;

$S_{14} = {}_p(714.2; 740.6) = \{719, 727, 733, 739\}$, #$S_{14} = 4$;
$S_{15} = {}_p(666.6; 689.5) = \{673, 677, 683\}$, #$S_{15} = 3$.

These Segments select a total of 1464 odd primes. We still have to choose 1556 – 1464 = 92 primes. The reader can view the remaining Segments and get an additional 60 primes. But closer to the origin, empty Segments appear and their number is growing rapidly. After the first hundred Segments, more and more time is spent on the analysis of empty ones. In this instance, the 345th Segment selects a prime factor 29, and the next prime 23 (descending) is selected by the 435th Segment. For small prime factors, it is impractical to iterate over Segments; it is faster to check directly prime numbers. □

Any prime is easy to check whether it falls into Chebyshev's Segment or into Black Hole. This analysis is useful for small primes as well as single-element Segments. Consider the corresponding algorithm.

**Algorithm 2.3.** Let $p$ be an odd prime. It is necessary to determine whether $p$ falls into any P-interval (2.2a) or not. If so, then $p \mid Cat(n)$. There are three steps in the algorithm.

*Step 1.* Calculate Chebyshev's Segment $S_u(n)$ whose lower boundary $(n+1)/u$ is as close as possible to $p$ from below. Easy to see, $u = \lceil (n+1)/p \rceil$.

*Step 2.* If $u \mid n+1$ (accordingly $p \mid n+1$ too), then $p \in H_u(n)$ and the algorithm is finished.

*Step 3.* It remains to compare $p$ and the upper boundary of the $u$-th Segment. If $p < 2n/(2u-1)$, then $p \in S_u(n)$. Otherwise $p \in H_{u-1}(n)$. □

Let's check the Algorithm 2.3 in the following example.

**Example 2.4.** For the 9,999-th Catalan number, let's check some primes selectively.

$p = 5$: $u = \lceil (9999+1)/5 \rceil = 2000$. Since $5 \mid 9999+1$, we get $5 \in H_{2000}$.
$p = 11$: $u = \lceil (9999+1)/11 \rceil = 910$. Since $11 > 2 \times 9999 /(2 \times 910 - 1) = 10.993$,
   we get $11 \in H_{909}$.
$p = 23$: $u = \lceil (9999+1)/23 \rceil = 435$. Since $23 < 2 \times 9999 /(2 \times 435 - 1) = 23.012$,
   we get $23 \in S_{435}$.
$p = 29$: $u = \lceil (9999+1)/29 \rceil = 345$. Since $29 < 2 \times 9999 /(2 \times 345 - 1) = 29.024$,
   we get $23 \in S_{435}$. □

**Factoring and … fishing**. Let the reader will forgive us for may be an inappropriate analogy. Let's compare the set of prime factors of the giant Catalan number with the vast sea. The prime factorization of the big Catalan numbers is reminiscent of fishing. Chebyshev's Segment like a fisherman's net catches a "flock" of big factors from the deep sea. But in shallow water near the shore small fish swims, the fishing net is useless there, and it is better to fish with a normal fishing rod. This is what the Algorithm 2.3 does; small factors are chosen one by one from the factor base of the Catalan number.

# 3 Legendre's layers

Chebyshev's Segments (2.2a) are able to select odd prime factors of *Cat*(*n*) over the P-interval $_p(2; 2n)$. The Segments do not intersect; if you combine all Segments then we will not get prime powers. But it is! For the prime number $p < \sqrt{2n}$, a power of 2 (squares) is possible and often occurs, a power of 3 (cubes) is possible and often occurs if $p < \sqrt[3]{2n}$, and so on. Obviously, Segments (2.2a) form a layer of the single (non-repeated) factors that are respectively less then 2*n*; we denote such a *SINGLE-layer* $L^{(1)}(n)$ and write $L^{(1)}(n) = \cup S_t^{(1)}(n)$, $t < (n+1)/2$. Next, we will talk about high layers: *SQUARE-layer* $L^{(2)}(n)$, *CUBE-layer* $L^{(3)}(n)$, and so on.

The power of *p* in *n*! is given by

(3.1) $$v_p(n!) = \sum_{k \geq 1} \lfloor n/p^k \rfloor, \quad n \in \mathbb{N},$$

where $\lfloor n/p^k \rfloor$ is the number of factors $p^k$ in $\{1, 2, \ldots, n\}$. In number theory, the last equality is known as *Legendre's formula* (some formulas for the p-adic valuation of the factorial see [Er19]).

Based on (1.1) and (3.1) we obtain

(3.2) $$v_p(Cat(n)) = \sum_{k \geq 1} \left( \lfloor 2n/p^k \rfloor - \lfloor n/p^k \rfloor - \lfloor (n+1)/p^k \rfloor \right).$$

In the canonical decomposition according to *k*, prime factors are distributed in layers, which we call *Legendre's layers*. In the case $k = 1$ we get SINGLE-layer $L^{(1)}(n)$. In each Legendre's layer, primes are selected by an own network of Chebyshev's Segments. Primes are not repeated in layers; a multiple prime factor of the Catalan number is selected from several layers.

Recall that for the *n*-th Catalan number, the elements of SINGLE-layer are selected from the P-interval $_p(2; 2n)$; and in $L^{(1)}(n)$ the first Segment $S_1^{(1)}(n) = {_p(n+1; 2n)}$ is obvious. Easy to see, the primes of SQUARE-layer are selected from the interval $_p(2; \sqrt{2n})$; and in $L^{(2)}(n)$ Segment $S_1^{(2)}(n) = {_p(\sqrt{n+1}; \sqrt{2n})}$ is also obvious. Note, if $p \in S_1^{(2)}(n)$ then $p \mid Cat(n)$, but not necessarily $p^2 \mid Cat(n)$. The first Segment is bordered by Black Hole $H_1^{(2)}(n) = {_p[\sqrt{2n/3}; \sqrt{n+1}]}$, which is preceded by the second Segment $S_2^{(2)}(n) = {_p(\sqrt{(n+1)/2}; \sqrt{2n/3})}$ and so on.

The general formula for Segments of the *n*-th SQUARE-layer is

$$S_t^{(2)}(n) = {_p(\sqrt{(n+1)/t}; \sqrt{2n/(2t-1)})}, \quad t \in \mathbb{N}.$$

The situation with the *n*-th CUBE-layer is similar. The elements in $L^{(3)}(n)$ are selected from the P-interval $_p(2; \sqrt[3]{2n})$, and for example, Segment $S_1^{(3)}(n) =$ $_p(\sqrt[3]{n+1}; \sqrt[3]{2n})$ is obvious. For the general case, we formulate the following theorem.

**Theorem 3.1.** *For the n-th Catalan number and the k-th Legendre's layer,*

$$S_t^{(k)}(n) = {}_p(\sqrt[k]{(n+1)/t}; \sqrt[k]{2n/(2t-1)}), \ t \in \mathbb{N}.$$

Obviously,

$$L^{(k)}(n) = \bigcup S_t^{(k)}(n), \ t \geq 1.$$

We use the obtained formulas for the decomposition of the Catalan number.

**Example 3.2.** Continue factoring the 9,999-th Catalan number. Let's select not empty Segments of the 9,999-th SQUARE-layer.

$S_1^{(2)} = {}_p(100; 141.4) = \{101, 103, \ldots, 139\}, \ \#S_1^{(2)} = 9;$

$S_2^{(2)} = {}_p(70.71; 81.64) = \{71, 73, 79\}, \ \#S_2^{(2)} = 3;$

$S_3^{(2)} = {}_p(57.73; 63.24) = \{59, 61\}, \ \#S_3^{(2)} = 2;$

$S_4^{(2)} = {}_p(50; 53.44) = \{51, 53\}, \ \#S_4^{(2)} = 2;$

$S_5^{(2)} = {}_p(44.72; 47.13) = \{47\}, \ \#S_5^{(2)} = 1;$

$S_6^{(2)} = {}_p(40.82; 42.63) = \{41\}, \ \#S_6^{(2)} = 1;$

$S_{12}^{(2)} = {}_p(28.86; 29.48) = \{29\}, \ \#S_{12}^{(2)} = 1;$

$S_{19}^{(2)} = {}_p(22.94; 23.24) = \{23\}, \ \#S_{19}^{(2)} = 1;$

$S_{28}^{(2)} = {}_p(18.89; 19.068) = \{19\}, \ \#S_{28}^{(2)} = 1;$

$S_{35}^{(2)} = {}_p(16.903; 17.024) = \{17\}, \ \#S_{35}^{(2)} = 1;$

$S_{83}^{(2)} = {}_p(10.976; 11.009) = \{11\}, \ \#S_{83}^{(2)} = 1.$

The elements of the 9,999-th SQUARE-layer (23 primes) are divided into 11 Segments. We calculated the first six Segments, but then we had to check each prime to accelerate (see below generalized Algorithm 3.3). The next CUBE-layer contains only three elements:

$S_1^{(3)} = {}_p(21.544; 27.143) = \{23\}, \ S_5^{(3)} = {}_p(12.599; 13.049) = \{13\},$ and
$S_8^{(3)} = {}_p(10.772; 11.006) = \{11\}.$

Below we present Segments of the remaining six layers.

$S_1^{(4)} = {}_p(10; 11.891) = \{11\};$

$S_1^{(5)} = {}_p(6.309; 7.247) = \{7\};$

$S_1^{(6)} = {}_p(4.641; 5.209) = \{5\}, \ S_{14}^{(6)} = {}_p(2.9898; 3.0079) = \{3\},$

$S_1^{(7)} = {}_p(3.727; 4.115) = \emptyset, \ S_5^{(7)} = {}_p(2.9619; 3.0068) = \{3\};$

$S_1^{(8)} = {}_p(3.162; 3.448) = \emptyset, \ S_2^{(8)} = {}_p(2.8998; 3.0059) = \{3\};$

$S_1^{(9)} = {}_p(2.7825; 3.0052) = \{3\}$.

As you can see, a prime number 3 has been repeated in layers 6, 7, 8, and 9. Therefore $v_3(Cat(9999)) = 4$. Note that in the 7th and 8th layers the first Segments are empty. □

Let's generalize Algorithm 2.3 to the case of any layer.

**Algorithm 3.3.** Let $p$ be an odd prime. It is necessary to determine whether $p$ falls into the $k$-th layer or not. If so, then $p \mid Cat(n)$.

*Step 1.* Calculate Segment $S_u^{(k)}(n)$ whose lower boundary $\sqrt[k]{(n+1)/u}$ is as close as possible to $p$ from below. Easy to see, $u = \lceil (n+1)/p^k \rceil$.

*Step 2.* If $p^k \mid n+1$ (accordingly $u \mid n+1$ too), then $p \in H_u^{(k)}(n)$ and the algorithm is finished.

*Step 3.* It remains to compare $p$ and the upper boundary of the $u$-th Segment. If $p < \sqrt[k]{2n/(2u-1)}$, then $p \in S_u^{(k)}(n)$. Otherwise $p \in H_{u-1}^{(k)}(n)$. □

# 4 Factorize by Kummer's theorem

In the process of factorization of the Catalan number, Chebyshev's Segments allow us to select sufficiently large groups of adjacent prime factors in each Legendre's layer. But it works well only at first. As a Segment number increases, the prime interval decreases dramatically. And to select the last 5-10% of prime factors, we have to spend almost all the time on the analysis of empty Segments.

Algorithm 3.3 directly checks a prime for getting into the Legendre's layer. But here, too, it is necessary to find the suitable Chebyshev's Segment and calculate its sometimes cumbersome boundaries. Consider an example with a giant number.

**Example 4.1.** In natural form, the Catalan number with index $10^8$ has a length of $6 \times 10^7$ decimal signs. Such a number is impossible to imagine, but it is not difficult to obtain a canonical decomposition. In SINGLE-layer, the smallest odd prime is 13, i.e. $13 \in S_u^{(1)}(10^8)$, where $u = \lceil (10^8+1)/13 \rceil = 7{,}692{,}308$. Let's calculate this Chebyshev's Segment:

$$S_{7,692,308}^{(1)}(10^8) = {}_p((10^8+1)/7{,}692{,}308;\ 2 \times 10^8/(2 \times 7{,}692{,}308 - 1))$$
$$= {}_p(12.9999996;\ 13.0000003).$$

A prime 13 also falls into SQUARE-layer, i.e. $13 \in L^{(2)}(10^8)$. The reader can calculate the corresponding Segment in the same way (in addition, it is necessary to work with square roots). □

The general view of Chebyshev's Segment allows us to modify well-known Kummer's theorem to factoring the Catalan numbers. Looking ahead, let's describe another data processing option from Example 4.1.

**Example 4.1a**. We calculate otherwise the divisibility of the $10^8$-th Catalan number by the prime 13 using modular arithmetic. Let's check two Legendre's layers.

SINGLE-layer: $13 \in L^{(1)}(10^8)$ as $10^8 \bmod 13 = 9 \in (13/2;\ 13-1) = (6.5;\ 12)$.

SQUARE-layer: $13 \in L^{(2)}(10^8)$ as $10^8 \bmod 13^2 = 165 \in (13^2/2;\ 13^2-1) = (64.5;\ 168)$.

Thus, $v_2(Cat(10^8)) \geq 2$. The logic of calculations is easy to see. In this case, we did not count Segments and did not extract square roots. □

**4.1. SINGLE-layer**. In Chebyshev's Segment (2.2a), the lower bound may be a prime, and then such a divisor doesn't fall into SINGLE-layer (but may be into another Legendre's layer). Thus, if an odd prime $p$ divides $n+1$, then $p \notin L^{(1)}(n)$. Accordingly, the following statement is true.

**Proposition 4.2.** *Let n be a positive integer and let an odd prime number p does not divide n+1. Then*

(i)   $\lfloor (n+1)/p \rfloor = \lfloor n/p \rfloor$;
(ii)  $n \bmod p < p-1$;
(iii) $\text{wt}_p(n+1) = \text{wt}_p(n) + 1$.

Recall $\text{wt}_p(n)$ denote the sum (*weight*) of the digits in the base-$p$ expansion of $n$.

Let's return to the sum (3.2), in which we are interested in the first term, i.e., the case $k = 1$. Taking into account Proposition 4.2(i), we obtain for an odd prime $p$

$$\lfloor 2n/p \rfloor - \lfloor n/p \rfloor - \lfloor (n+1)/p \rfloor = \lfloor 2n/p \rfloor - 2\lfloor n/p \rfloor,\ p \nmid n+1.$$

Let $\{x\} = x - \lfloor x \rfloor$, the fractional part of a real number $x$. Then (see [Po15], p. 2)

(4.1) $\quad \lfloor 2n/p \rfloor - 2\lfloor n/p \rfloor = 2\{n/p\} - \{2n/p\} = \begin{cases} 1, & \text{if } \{n/p\} > \tfrac{1}{2}, \\ 0, & \text{if } \{n/p\} < \tfrac{1}{2}. \end{cases}$

Obviously, for an integer $n$ and an odd $p$, $\{n/p\} \neq \tfrac{1}{2}$. According to (4.1), a prime number $p$ falls into SINGLE-layer if we get the carry by doubling $\{n/p\}$ under the condition that $n \bmod p < p-1$ (i.e. $p \nmid n+1$). In fact, we have echoes of Kummer's theorem from 1852 [Po15].

Let's issue the received result in the form of the theorem.

**Theorem 4.3.** *For the n-th Catalan number, an odd prime p falls into SINGLE-layer if and only if*

$$p/2 < n \bmod p < p-1 \quad \text{or} \quad n \bmod p \in (p/2, p-1).$$

As a result, SINGLE-layer is easily filled with prime factors. For example, a prime 5 falls into $L^{(1)}(n)$ if and only if $n \bmod 5 = 3 \in (5/2;\ 5-1) = (2\tfrac{1}{2},\ 4)$ or $n = 5i + 3$, $i \geq 0$. For a prime 7, we get $n = 7i + (4$ or $5)$. What about a prime 3? It's simple, there are no integers in the interval $(3/2;\ 3-1) = (1½,\ 2)$. Hence

**Corollary 4.4.** *For the n-th Catalan number,*

(i) $\quad 3 \notin L^{(1)}(n)$;

(ii) $\quad 5 \in L^{(1)}(n)$ *if and only if* $n = 5i + 3$, $i \geq 0$;

(iii) $\quad 7 \in L^{(1)}(n)$ *if and only if* $n = 7i + (4$ or $5)$, $i \geq 0$.

**4.2. The *k*-th Legendre's layer.** Theorem 3.1 gives a general view of Chebyshev's Segment for the *k*-th Legendre's layer. The lower boundary of such Segment can be an integer and, accordingly, a prime number. For the *n*-th Catalan number, let a prime $p = \sqrt[k]{(n+1)/t} \notin L^{(k)}(n)$. Then

$$p^k = (n+1)/t, \quad t = (n+1)/p^k, \text{ and hence } p^k \mid n+1.$$

Now, let's generalize the previous statement.

**Proposition 4.5** (generalization of Proposition 4.2). *Let n be a positive integer, let p be an odd prime number and let $p^k \nmid n+1$, $k \in \mathbb{N}$. Then*

(i) $\quad \lfloor (n+1)/p^k \rfloor = \lfloor n/p^k \rfloor$;

(ii) $\quad n \bmod p^k < p^k - 1$;

(iii) $\quad \mathrm{wt}_{p^k}(n+1) = \mathrm{wt}_{p^k}(n) + 1$.

Next, we will follow Carl Pomerance's analysis [Po15] for binomial coefficients. Let's use the Proposition 4.5(i) and simplify the *k*-th summand in (3.2).

$$\lfloor 2n/p^k \rfloor - \lfloor n/p^k \rfloor - \lfloor (n+1)/p^k \rfloor = \lfloor 2n/p^k \rfloor - 2\lfloor n/p^k \rfloor, \quad p \in Lev_k(n).$$

Based on the equality $n/p^k = \lfloor n/p^k \rfloor + \{n/p^k\}$, let's move to the fractional parts

(4.2) $\qquad 2\{n/p^k\} - \{2n/p^k\} = \begin{cases} 1, & \text{if } \{n/p^k\} > ½, \\ 0, & \text{if } \{n/p^k\} < ½. \end{cases}$

Obviously, $\{n/p^k\} \neq ½$, since $p$ is odd. Based on (4.2) and Proposition 4.5 we obtain a generalized theorem.

**Theorem 4.6** (generalization of Theorem 4.3). *For the n-th Catalan number, an odd prime p falls into the k-th Legendre's layer if and only if*

$$½p^k < n \bmod p^k < p^k - 1 \quad \text{or} \quad n \bmod p^k \in (½p^k, p^k - 1).$$

The resulting Theorem 4.6 can be called a modification of Kummer's theorem for the *k*-th Legendre's layer.

**Example 4.7**. Continue factoring the Catalan number with index $n=10^8$. Let's find all the prime factors of the 7th Legendre's layer, whose elements are selected from the following P-interval:

$$_p(2; \sqrt[7]{2n}) = {}_p(2; 15.34) = \{3, 5, 7, 11, 13\}.$$

We need to check 5 primes.

Prime 3:   $10^8 \bmod 3^7 = 1{,}612 \in (3^7/2; 3^7-1) = (1{,}093.5; 2{,}186);$  so $3 \in L^{(7)}$.

Prime 5:   $10^8 \bmod 5^7 = 0 \notin (5^7/2; 5^7-1) = (39{,}062.5; 78{,}124);$  so $5 \notin L^{(7)}$.

Prime 7:   $10^8 \bmod 7^7 = 351{,}297 \notin (7^7/2; 7^7-1) = (411{,}771.5; 823{,}542);$  so $7 \notin L^{(7)}$.

Prime 11: $10^8 \bmod 11^7 = 2{,}564{,}145 \notin (9{,}743{,}585.5; 19{,}487{,}170);$  so $11 \notin L^{(7)}$.

Prime 13: $10^8 \bmod 13^7 = 37{,}251{,}483 \in (31{,}374{,}258.5; 62{,}748{,}516);$  so $13 \in L^{(7)}$.

Thus, there are only two elements in the 7th layer: $L^{(7)}(100{,}000{,}000) = \{3, 13\}$.   □

## 5   Online software service.

The reader can independently factorize any Catalan number and check certain results using this program. Layer-by-layer factorization is performed with segmentation of the SINGLE-layer without issuing a naturalized Catalan number. Often control calculations allow us to conduct independent research. The program code in the form of HTML-file is available to everyone. The participation of those wishing to improve the software service is welcome.

**Acknowledgements.** I wish to thank Igor Pak (Mathematics Department at UCLA, USA) for the detailed and comprehensive historical review of Catalan-like numbers, and thanks for a selection of recent works on the subject. The author would like to thank and express his sincere gratitude to Bruce Sagan (Michigan State University, USA) for helpful comments.

Gzhel State University, Moscow, 140155, Russia

http://www.en.art-gzhel.ru/